\documentclass{amsart}
\usepackage{graphicx}
\usepackage{latexsym}
\usepackage{hyperref}
\usepackage{amsfonts}
\usepackage[all]{xy}
\usepackage{amssymb, mathrsfs, amsfonts, amsmath}
\usepackage{amsbsy}
\usepackage{amsfonts}
\newtheorem{remark}{Remark}

\newtheorem{theorem}{Theorem}

\numberwithin{equation}{section}

\newcommand{\disp}{\displaystyle}

\def\Om{{\Omega}}

\def\re{{\mathbb{R}}}

\def\Om{{\Omega}}

\def\jnt{\disp\int}
\def\jjnt{\jnt\!\!\!\jnt}

\newcommand{\Fin}{\hfill$\Box$}
\numberwithin{equation}{section}

\begin{document}
\title[hierarchic control for the wave equation]{Remarks on hierarchic control for the wave equation in a non cylindrical domain}
\author{Isa\'{\i}as  Pereira de Jesus$^{1}$}
\address{$^{1}$ Departamento de Matem\'{a}tica-Universidade Federal do Piau\'i, 64049-550-Teresina /PI, Brazil}
\email{isaias@ufpi.edu.br} 
\date{April 15, 2014}


\begin{abstract}
In this paper we establish hierarchic control for the wave equation in  a non
cylindrical domain $\widehat{Q}$  of  $\mathbb{R}^{n + 1}$. We assume that we
can act in the dynamic of the system by a hierarchy of controls. According to
the formulation given by   H. Von Stackelberg \cite{S}, there are local controls,
called followers and global controls, called leaders. In fact, one considers
situations where there are two cost (objective) functions. One possible way is to
cut the control into two parts, one being thought of as "the leader" and the other
one as "the follower". This situation is studied in the paper, with one of the cost
functions being of the controllability type. Existence and uniqueness is proven.
The optimality system is given in the paper.

\end{abstract}

\maketitle
\section{Introduction}\label{introd}
Optimization problems appear quite often in a number of problems from engineering sciences and economics. Many mathematical models are formulated in terms of optimization problems involving a unique objective, minimize
 cost or maximize profit, etc... In more realistic and involved situations, several (in general conflicting) objectives must be considered. In general for the classical mono-objective control problem, a functional adding the objectives
 of the problem is defined and an unique control is used. When it is not clear how to average the different objectives or when several controls are used,
it becomes essential to study multi-objective problems. The different notions of
equilibrium  for multi-objective problems were introduced in economy and games
theory (see \cite{N, P, S}).

 In the context of game theory, the Stackelberg strategy is normally applied within games where some players are in a better position than others. The dominant players are called
 the leaders and the sub-dominant players the followers. One situation where this concept is very convincing is when the players choose their strategies one after another and the player
 who does the first move has an advantage.
This interpretation however is not so convincing in continuous time where players
choose theories strategies at each time, but one could think of that one player
has a larger information structure. However from a mathematical point the of
view, the  concept is so successful because it leads to results.

This paper was inspired by the ideas contained in the work of J.-L. Lions
\cite{L1}, where we  investigate a  similar question of hierarchic control
employing the Stackelberg strategy in the case of time dependent domains.

 We also mention some previous works related to our theme. In the articles by
Lions  \cite {L10, L11, L12} the author gives some results on Pareto and
Stackelberg equilibrium. D\'iaz and Lions \cite {DL} studied the existence of
Stackelberg-Nash equilibrium. Similar problems for the Stokes equations are
considered in Gonz\'alez, Lopes and Medar \cite {GO}. On the other hand,
Glowinski, Ramos and Periaux \cite {RA1, RA2} studied the Nash equilibrium
from a theoretical and numerical point of view, first for linear parabolic differential
equations \cite {RA1} and afterwards for the Burgers equation \cite {RA2}. In
Limaco, Clark and Medeiros \cite {LI}, the authors present the Stackelberg-Nash
equilibrium in the context of linear heat equation in non cylindrical domains. In
Russell \cite {Ru}, the author present the boundary value control of the
higher-dimensional wave equation.  We can mention also Jesus\cite{Je},  in
which the author studied the hierarchic control for the wave equation in moving
domains.

We will obtain the following three main results: the existence and  uniqueness of
Nash equilibrium, the approximate controllability  with respect to the leader
control and the optimality system for the leader control.

The rest of the paper is organized as follows. In Section \ref{sec2} the formulation
of  problem is established. The existence and uniqueness of the Nash
equilibrium is deduced in Section \ref{sec3}. In Section \ref{sec4} we study  the
approximate controllability with respect to the leader control. Finally, in the
Section \ref{sec5} we present  the optimality system for the leader control.

\section{Problem formulation}\label{sec2}

As in \cite{Mi}, let $\Om$ be an open  bounded set of $\mathbb{R}^n$ with
boundary $\Gamma$ of class $C^2$, which without  loss of gererality, can be
assumed containing the origin of  $\mathbb{R}^n$, and  $k : [0,\infty) \mapsto (0,
\infty)$ a continuously differentiable function. Let us consider the subsets
$\Om_t$ of $\mathbb{R}^n$ given by
$$ \Om_t = \{x \in \mathbb{R}^n; x = k(t)y, y \in \Om\}, 0 \leq t \leq T < \infty,$$
whose boundaries are denoted by $\Gamma_t$ and $\widehat{Q}$ the non
cylindrical domain of $\mathbb{R}^{n + 1}$,
$$\disp \widehat{Q}= \bigcup_{0 \leq t \leq T} \{ \Om_t \times \{t\}\}$$
with lateral boundary  defined by $\disp \widehat{\Sigma} =\widehat{\Sigma}_0
\cup \widehat{\Sigma}_0^*$, where
$$\disp \widehat{\Sigma} = \bigcup_{0 \leq t \leq T} \{ \Gamma_t \times \{t\}\}.$$

Thus we consider the mixed problem
\begin{equation} \label{eq1.3}
\left|
\begin{array}{l}
\displaystyle u'' - u_{xx} = 0 \ \ \mbox{ in } \ \ \widehat{Q}\\[11pt]
\disp u(x,t) = \left\{
\begin{array}{l}
\widetilde{w} \ \ \mbox{ on } \ \ \widehat{\Sigma}_{0}\\[11pt]
0 \ \ \mbox{ on } \ \ \widehat{\Sigma}_0^*
\end{array}
\right.\\[11pt]
\disp u(x,0) = u_0(x), \;\; u'(x,0) = u_1(x) \ \mbox{ in }\; \Om_0.
\end{array}
\right.
\end{equation}

By $\disp u''=u''(x,t)$ we represent the derivative $\disp \frac{\partial^2 u}{\partial
t^2}$ and $\disp u_{xx}=u_{xx}(x,t)$ the Laplace operator $\disp \frac{\partial^2
u}{\partial x^2}$.

In this article, motivated by the arguments contained in the work of J.-L. Lions
\cite{L1}, we investigate a  similar question of hierarchic control for the equation
(\ref{eq1.3}), employing the Stackelberg strategy in the case of time
 dependent domains.

Following the work of J.-L. Lions \cite{L1}, we divide $\disp \widehat{\Sigma}_0$
in two disjoint  parts
\begin{equation}\label{decomp0}
\disp \widehat{\Sigma}_0=\disp \widehat{\Sigma}_1 \cup \disp \widehat{\Sigma}_2,
\end{equation}
 and consider
\begin{equation} \label{decomp}
\disp \widetilde{w}=\{\widetilde{w}_1, \widetilde{w}_2\}, \;\; \widetilde{w}_i=\mbox{control function in } \; L^2(\widehat{\Sigma}_i), \;i=1,2.
\end{equation}

Thus, we observe  that the system (\ref{eq1.3}) can be rewritten as follows:
\begin{equation} \label{eq1.3.1}
\left|
\begin{array}{l}
\displaystyle u''- u_{xx} = 0 \ \ \mbox{ in } \ \ \widehat{Q}\\[5pt]
\disp u(x,t) = \left\{
\begin{array}{l}
\widetilde{w_1} \ \ \mbox{ on } \ \ \widehat{\Sigma}_{1}\\
\widetilde{w_2} \ \ \mbox{ on } \ \ \widehat{\Sigma}_{2}\\
0 \ \ \mbox{ on } \ \ \widehat{\Sigma}\backslash \widehat{\Sigma}_0
\end{array}
\right.\\[5pt]
\disp u(x,0) = u_0(x), \;\; u'(x,0) = u_1(x) \ \mbox{ in }\; \Om_0.
\end{array}
\right.
\end{equation}
In the decomposition (\ref{decomp0}), (\ref{decomp}) we establish a hierarchy.
We think of $\widetilde{w}_1$ as being the "main" control, the leader  and we
think of $\widetilde{w}_2$ as the follower,  in Stackelberg terminology.

Associated to the solution $u=u(x,t)$ of (\ref{eq1.3.1}), we will consider the
(secondary) functional
\begin{equation}\label{sfn}
\disp \widetilde{J}_2(\widetilde{w}_1, \widetilde{w}_2) =
\displaystyle\frac{1}{2} \displaystyle\jjnt_{\widehat{Q}} \left(u(\widetilde{w}_1, \widetilde{w}_2)-\widetilde{u}_2\right)^2 dx dt +
\displaystyle\frac{\widetilde{\sigma}}{2} \int_{\widehat{\Sigma}_2} \widetilde{w}_2^2\;d\widehat{\Sigma},
\end{equation}
and the (main) functional
\begin{equation}\label{mfn}
\disp \widetilde{J}(\widetilde{w}_1) =\frac{1}{2}\int_{\widehat{\Sigma}_1} \widetilde{w}_1^2\;d\widehat{\Sigma},
\end{equation}
where $\widetilde{\sigma}>0$ is constant and $\widetilde{u}_2$ is a given
function in $L^2(\widehat{Q}).$

\begin{remark}\label{bdfnc} From the regularity and uniqueness of the solution  of (\ref{eq1.3.1}) (see Remark \ref{rsol} ) the cost functionals $\disp \widetilde{J}_2$ and $\disp \widetilde{J}$ are well defined.
\end{remark}
The control problem that we will consider is as follows: the follower $\disp
\widetilde{w}_2$ assume that  the leader $\disp \widetilde{w}_1$ have made  a
choice. Then, it tries  to find an equilibrium of the cost  $\widetilde{J_2}$ , that is,
it  looks for  a control $\disp \widetilde{w}_2=\mathfrak{F}(\widetilde{w}_1)$
(depending on $\disp \widetilde{w}_1$), satisfying:
\begin{equation}\label{n1.1}
\disp \widetilde{J}_2(\widetilde{w}_1, \widetilde{w}_2) \leq \widetilde{J}_2(\widetilde{w}_1,  \widehat{w}_2), \;\;\; \forall\; \widehat{w}_2  \in L^2(\widehat{\Sigma}_2).
\end{equation}

The control $ \widetilde{w}_2 $, solution of the (\ref{n1.1}), is called Nash
equilibrium for the cost $\widetilde{J}_2 $ and it depends on $ \widetilde{w}_1$
(cf. Aubin \cite{A}).

\begin{remark}\label{r2} In another way, if the leader $ \widetilde{w}_1$ makes a choice, then the follower $ \widetilde{w}_2$ makes also a choice, depending on  $\widetilde{w}_1$, which becomes minimum it cost
 $\widetilde{J}_2$, that is,
\begin{equation}\label{son}
\disp \widetilde{J}_2(\widetilde{w}_1, \widetilde{w}_2)= \inf_{\widehat{w}_2 \in L^2(\widehat{\Sigma}_2)} \widetilde{J}_2(\widetilde{w}_1, \widehat{w}_2).
\end{equation}

This is equivalent to the  (\ref{n1.1}). This process is called Stackelberg-Nash
strategy, see D\'iaz and Lions \cite{DL}.
\end{remark}

After this, we consider the state $\disp u\left(\widetilde{w}_1,
\mathfrak{F}(\widetilde{w}_1)\right)$ given by the solution of
\begin{equation} \label{eq1.3.1F}
\left|
\begin{array}{l}
\displaystyle u''- u_{xx} = 0 \ \ \mbox{ in } \ \ \widehat{Q}\\[3pt]
\disp u(x,t) = \left\{
\begin{array}{l}
\widetilde{w_1} \ \ \mbox{ on } \ \ \widehat{\Sigma}_{1}\\
\mathfrak{F}(\widetilde{w}_1) \ \ \mbox{ on } \ \ \widehat{\Sigma}_{2}\\
0 \ \ \mbox{ on } \ \ \widehat{\Sigma}\backslash \widehat{\Sigma}_0
\end{array}
\right.\\[3pt]
\disp u(x,0) = u_0(x), \;\; u'(x,0) = u_1(x) \ \mbox{ in }\; \Om_0.
\end{array}
\right.
\end{equation}

We will look for any optimal control $\disp \widetilde{w_1}$ such that
\begin{equation}\label{ocn}
\disp \widetilde{J}(\widetilde{w_1}, \mathfrak{F}(\widetilde{w}_1))= \inf_{\overline{w}_1 \in L^2(\widehat{\Sigma}_1)} \widetilde{J}(\overline{w}_1, \mathfrak{F}(\widetilde{w}_1)),
\end{equation}
subject to the restriction of the approximate controllability type
\begin{equation}\label{apcn}
\begin{array}{l}
\disp\left( u(x, T; \widetilde{w_1}, \mathfrak{F}(\widetilde{w}_1)), u'(x, T; \widetilde{w_1}, \mathfrak{F}(\widetilde{w}_1))\right)\in B_{L^2(\Omega_t)}(u^0,\rho_0) \times B_{H^{-1}(\Omega_t)}(u^1,\rho_1),
\end{array}
\end{equation}
 where $\disp B_X(C, r)$ denotes the ball in $X$ with center $C$ and ratio $r$.

To explain this optimal problem, we are going to consider the following
sub-problems:

\textbf{$\bullet$ Problem 1} Fixed any leader control $\widetilde{w}_1$,  find the
follower control $\disp \widetilde{w}_2=\mathfrak{F}(\widetilde{w}_1)$ (depending
on $\disp \widetilde{w}_1$) and the associated state $u$, solution of
(\ref{eq1.3.1}) satisfying the condition (\ref{son}) (Nash equilibrium)  related to
$\widetilde{J}_2$, defined in (\ref{sfn}).

\textbf{$\bullet$ Problem 2} Assuming that the existence of the Nash equilibrium
$\disp \widetilde{w}_2$ was proved, then when $\widetilde{w_1}$ varies in
$L^2(\Omega_t)$,  prove that the solutions $\disp\left(u(x, t; \widetilde{w}_1,
\widetilde{w}_2), u'(x, t; \widetilde{w}_1, \widetilde{w}_2)\right)$ of the state
equation (\ref{eq1.3.1}), evaluated at $t = T$, that is, $\disp\left(u(x, T;
\widetilde{w}_1, \widetilde{w}_2), u'(x, T; \widetilde{w}_1, \widetilde{w}_2)\right)$,
generate a dense subset of $\disp L^2(\Omega_t) \times H^{-1}(\Omega_t)$.

\begin{remark}\label{r1} By the linearity of the system (\ref{eq1.3.1F}), without loss of generality we may assume that $u_0=0=u_1$.\\
\end{remark}
Following the work of J.-L. Lions \cite{L1}, we divide $\disp {\Sigma}_0$ in two
disjoint parts
\begin{equation}\label{decomp0.1}
\disp {\Sigma}_0=\disp {\Sigma}_1 \cup \disp {\Sigma}_2,
\end{equation}
and consider
\begin{equation} \label{decomp.2}
\disp w = \{w_1,w_2\}, \;\; {w}_i=\mbox{control function in } \; L^2({\Sigma}_i), \;i=1,2.
\end{equation}
We can also write
\begin{equation} \label{decomp 2.A}
\disp w= w _1+ w_2, \; \mbox{ with } \; \disp {\Sigma}_0=\disp {\Sigma}_1 = \disp {\Sigma}_2.
\end{equation}

Note that when $\disp (x,t)$ varies in $\disp \widehat{Q}$ the point $\disp (y,t)$,
with $\disp y=\frac{x}{k(t)}$, varies in $\disp Q=\Om \times (0,T)$. Then, the
application $$\disp \zeta:\widehat{Q} \to Q, \;\;\;\zeta(x,t)=(y,t)$$ is of class
$\disp C^2$ and the inverse $\disp \zeta^{-1}$ is also of class $\disp C^2$.
Therefore the change of variables $\disp u(x,t)=v(y,t)$, transforms the
initial-boundary value problem (\ref{eq1.3.1}) into the equivalent system
\begin{equation} \label{eq1.14}
\left|
\begin{array}{l}
\disp v'' + Lv = 0  \ \ \mbox{ in } \ \ Q\\ [13pt]
\disp v(y,t) = \left\{
\begin{array}{l}
w_1 \ \ \mbox{on} \ \ \Sigma_1\\
w_2 \ \ \mbox{on} \ \ \Sigma_2\\
0 \ \ \mbox{on} \ \ \Sigma\backslash\Sigma_0
\end{array}
\right. \\[5pt]
\disp v(y,0) = v_0(y), \; v'(y,0) = v_1(y), \;\; y \in \Om,
\end{array}
\right.
\end{equation}
where
\begin{equation*}
\left|
\begin{array}{l}
\disp Lv = - [a_{ij}(y,t)v_{y_j}]_{y_i} + b_i(y,t)v'_{y_i} + c_i(y,t)v_{y_i} \\[10pt]
\disp a_{ij}(y,t)= (\delta_{ij} - k'^2y_iy_j)k^{-2}\\[10pt]
\disp b_i(y,t) = -2k'k^{-1}y_i\\[10pt]
\disp c_i(y,t) = [(1 - n)k'^{2} - k''k]k^{-2}y_i \\[10pt]
\disp \Sigma= {\Sigma}_0 \cup \Sigma^*_0\\[10pt]\disp
\disp {\Sigma}_0 = \Gamma_{0} \times (0,T)\\[10pt]
\disp {\Sigma}^*_0 = \Sigma\backslash\Sigma_0\\[10pt]
\end{array}
\right.
\end{equation*}

We consider the coefficients of the operator $L$  satisfying  the following
conditions :
\begin{itemize}
\item[(H1)] $a_{ij}(y,t)$ are symmetric and uniformly coercive on $Q$;
\item[(H2)] $a_{ij}(y,t) \in C^{1}(\overline{Q})$,   $a''_{ij}(y,t)  \in  L^{\infty}(Q)$;
\item[(H3)] $b_i(y,t), c_i(y,t)  \in  W^{1,\infty}(0, T ; L^{\infty}(\Om))$, $(b_i(y,t)_{y_i} \in L^{\infty}(Q)$,
\end{itemize}

and concerning the function $k$,
\begin{itemize}
\item[(H4)] $k \in  W^{3,\infty}_{loc}(0, \infty).$
\end{itemize}

In this way, it is enough to investigate the control problem for the equivalent
problem (\ref{eq1.14}).

$\bullet$ \textbf{Cost functionals in the cylinder $Q$.}  From the diffeomorphism
$\disp \zeta$ which transforms $\widehat{Q}$ in $Q$, we transform the cost
functionals $\disp \widetilde{J}_2, \disp \widetilde{J}$ in the cost functionals
$J_2, J$ defined by
\begin{equation} \label{eq3.9}
\disp J_2(w_1,w_2) = \frac{1}{2}\,\int_{0}^{T}\int_{\Omega}k^n(t)[v(w_1,w_2) - v_2(y,t)]^2dy\,dt + \frac{\sigma}{2}\,\int_{\Sigma_2}w_{2}^{2}\,d\Sigma
\end{equation}
and
\begin{equation} \label{eq3.7}
\disp J(w_1) = \frac{1}{2}\,\int_{\Sigma_1}w_1^{2}\,d\Sigma,
\end{equation}
where $\sigma > 0$ is constant and $v_2(y,t)$ is a given function in $L^2(\Om
\times (0,T)).$
\begin{remark}\label{rsol} By using similar technique as used in \cite{Mi}, we can prove the  following:
For each $v_0 \in H_0^1(\Om)$ and $v_1 \in L^2(\Om)$, and $\disp w_i \in
L^2(\widehat{\Sigma}_i), \;i=1,2$ there exists exactly one  solution $v$ to
(\ref{eq1.14}), with  $\disp v \in C\big( [0,T] ; H_0^1(\Om)) \cap C^{1}\big( [0,T] ;
L^2(\Om))$, see \cite{Mi}. Thus, in particular, the cost functionals $\disp {J}_2$
and $\disp {J}$ are well defined.

By using the diffeomorphism   $\disp \zeta^{-1}(y,t)=(x,t)$, from $Q$ to
$\widehat{Q}$, we obtain an unique global weak solution $\disp u$ for the
problem (\ref{eq1.3.1}) with the regularity $\disp u \in C\big( [0,T] ; H_0^1(\Om_t))
\cap C^{1}\big( [0,T] ; L^2(\Om_t))$.
\end{remark}

Associated to functionals $\disp J_2$ and $\disp J$ defined above, we will
consider the following  sub-problems:

\textbf{$\bullet$ Problem 3} Fixed any leader control ${w}_1$,  find the  follower
control $\disp {w}_2$ (depending on $\disp {w}_1$) and the associated state $v$
solution of (\ref{eq1.14}) satisfying (Nash equilibrium)
 \begin{equation}\label{soncil}
\disp {J}_2({w}_1, {w}_2)= \inf_{\widehat{w}_2 \in L^2({\Sigma}_2)} {J}_2({w}_1, \widehat{w}_2),
\end{equation}
related to ${J}_2$ defined in (\ref{eq3.9}).

\textbf{$\bullet$ Problem 4} Assuming that the existence of the Nash equilibrium
$\disp {w}_2$ was proved, then when ${w_1}$ varies in $L^2(\Omega)$,  prove
that the solutions $\disp\left(v(y, t; {w_1}, {w}_2), v'(x, t; {w_1}, {w}_2)\right)$ of
the state equation (\ref{eq1.14}), evaluated at $t = T$, that is, $\disp\left(v(y, T;
{w_1}, {w}_2), v'(y, T; {w_1}, {w}_2)\right)$, generate a dense subset of $\disp
L^2(\Om) \times H^{-1}(\Om)$.

 \section{Nash equilibrium}\label{sec3}
In this section, fixed any leader control $\disp \ w_1 \in L^2(\Sigma_1)$ we want
to determine the existence and uniqueness of solutions for the problem
\begin{equation} \label{eq3.10}
\begin{array}{l}
\displaystyle\inf_{w_2 \in L^2(\Sigma_2)}J_2(w_1,w_2),
\end{array}
\end{equation}
and then to obtain a characterization of this solution in terms of an adjoint
system.

In fact, this is a classical type problem in the control of distributed systems (cf.
J. L. Lions \cite{L3}). It admits an unique solution
\begin{equation} \label{eq3.11}
\disp w_2 = \mathfrak{F}(w_1).
\end{equation}
The Euler equation for (\ref{eq3.10}) is given by
\begin{equation} \label{eq3.21}
\int_{0}^{T}\int_{\Omega}{k^n(t)}(v-v_2)\widehat{v}dy\,dt + \sigma\int_{\Sigma_2}w_2\widehat{w}_2d\Sigma = 0, \;\;\forall\, \widehat{w}_2 \in L^2(\Sigma_2)
\end{equation}
where $\widehat{v}$ is solution of the following system
\begin{equation} \label{eq3.22}
\left|
\begin{array}{l}
\disp \widehat v'' +L\widehat{v} = 0 \ \ \mbox{ in } \ \ Q\\ [7pt]
\disp \widehat{v} = \left\{
\begin{array}{l}
0 \ \ \mbox{ on } \ \ \Sigma_1 \\
\widehat{w}_2 \ \ \mbox{ on } \ \ \Sigma_2 \\
0 \ \ \mbox{ on } \ \ \Sigma \backslash \left(\Sigma_1 \cup \Sigma_2\right)
\end{array}
\right. \\[13pt]
\disp \widehat{v}(y,0) = 0, \; \widehat{v'}(y,0) = 0, \;\; y \in \Om.
\end{array}
\right.
\end{equation}
In order to express (\ref{eq3.21}) in a convenient form, we introduce the adjoint
state defined by
\begin{equation}\label{sac}
\left|
\begin{array}{l}
p'' + L^{\ast}\,p = k^n(t)\left(v - v_2\right) \ \mbox{ in } \ \ Q, \\[5pt]\disp
p(T) = p'(T) = 0, \;\; y \in \Om, \\[5pt]\disp
p = 0 \ \mbox{ on } \ \Sigma,
\end{array}
\right.
\end{equation}
where
\begin{equation*}
\begin{array}{l}
\disp L^{\ast}\,p =  - \left[ (\delta_{ij} - k'^2 y_i y_j)k^{-2}p_{y_j}\right]_{y_i} - 2k'k^{-1}y_ip'_{y_i} - 2nk'k^{-1}p'\\[9pt] \disp
+ \left[(n + 1)k'^2 - k'' k\right] k^{-2}y_ip_{y_i} +  \left[n(n + 1)k'^2 - nk''k \right] k^{-2}p
\end{array}
\end{equation*}
is the formal adjoint of the  operator $\disp L$.

If we multiply (\ref{sac}) by $\widehat{v}$ and if we integrate by parts, we find
\begin{equation} \label{eq3.33}
\int_{0}^{T}\int_\Om k^n(t)(v - v_2)\widehat{v}\,dy\,dt + \int_{\Sigma_2} \frac{1}{k^2(t)}\,p_\nu\,\widehat{w}_2\,d\Sigma = 0,
\end{equation}
so that (\ref{eq3.21}) becomes
\begin{equation}\label{ci}
\disp p_\nu= \sigma k^2(t)\,w_2 \ \ \mbox{ on } \ \ \Sigma_2.
\end{equation}
We can summarize as follows:
\begin{theorem}\label{teN} For each $\disp w_1 \in L^2(\Sigma_1)$ there exists an unique Nash equilibrium $\disp w_2$ in the sense of (\ref{soncil}). Moreover, the follower $\disp w_2$ is given by
\begin{equation}\label{cseg}
\disp \disp w_2 = \mathfrak{F}(w_1)=\frac{1}{\sigma k^2(t)}\,\;p_\nu\;\;\mbox{ on }\;\;\Sigma_2,
\end{equation}
where $\disp \{ v,p \}$ is the unique solution of (the optimality system)
\begin{equation} \label{eq3.37}
\left|
\begin{array}{l}
\disp v'' + Lv = 0 \ \mbox{ in } \;\; Q \\
\disp p'' + L^{\ast}\,p = k^n(t)\left(v - v_2\right) \ \mbox{ in } \;\; Q\\[5pt]
\disp v = \left\{
\begin{array}{l}
w_1 \ \mbox{ on } \ \Sigma_1\\[5pt]
\disp \frac{1}{\sigma k^2(t)}\;\,p_\nu \ \mbox{ on } \ \Sigma_2\\[5pt]
0 \ \mbox{ on } \ \Sigma \backslash \Sigma_0
\end{array}
\right.\\[7pt]
p = 0 \ \mbox{ on } \ \Sigma \\[5pt]
v(0) = v'(0) = 0, \;\; y \in \Om\\[5pt]
p(T) = p'(T) = 0, \;\; y \in \Om.
\end{array}
\right.
\end{equation}
Of course $\disp \{ v,p \}$ depends on $w_1$:
\begin{equation}\label{cdep}
\disp \{ v,p \} = \{ v(w_1),p(w_1) \}.
\end{equation}
\end{theorem}

\section{On the approximate controllability}\label{sec4}
Since we have proved the existence, uniqueness and characterization of the
follower $\disp w_2$, the leader $\disp w_1$  now wants  that the solutions $v$
and $v'$, evaluated at time $t=T$, to  be as close as possible to $\disp(v^0,
v^1)$. This will be possible if the system (\ref {eq3.37})  is approximate
controllable. We are looking for
\begin{equation} \label{inf1}
\begin{array}{l}
\displaystyle\inf\, \frac{1}{2\,}\,\int_{\Sigma_1} w_{1}^{2}\,d\Sigma
\end{array}
\end{equation}
where $\disp w_1$ is subject to
\begin{equation} \label{subj1}
\begin{array}{l}
\disp\left(v(T;{w_1}), v'(T; {w_1})\right)\in B_{L^2(\Om)}(v^0,\rho_0) \times B_{H^{-1}(\Om)}(v^1,\rho_1),
\end{array}
\end{equation}
assuming that  $w_1$ exist, $\rho_0$ and  $\rho_1$  positive numbers arbitrarily
small  and $\{v^0, v^1\} \in L^2(\Om) \times H^{-1}(\Om)$.

To study (\ref{inf1}), we need the hypotheses
\begin{equation} \label{hT}
\disp T>2d(\Om,\Gamma_0), \;\; \mbox{ where }\;\; \disp d(\Om, \Gamma_0) = \sup_{x \in \Om} d(x,\Gamma_0)
\end{equation}

Now as in the case (\ref{decomp 2.A}) and using  Holmgren's Uniqueness
Theorem (cf. \cite{LH})  and for the explicit use of it made here (cf. \cite {Mi}), the
following approximate controllability result holds:
\begin{theorem}\label{AC} Assume that (\ref{hT}) hold. Let us consider $\disp w_1 \in L^2(\Sigma_1)$ and $\disp w_2$ a Nash equilibrium in the sense (\ref{soncil}).
Then the functions $\disp \left(v(T), v'(T)\right)=\left(v(., T, {w_1}, w_2), v'(., T,
{w_1}, w_2)\right)$, where $\disp v$ solves the system (\ref{eq1.14}), generate a
dense subset of $\disp L^2(\Om)\times H^{-1}(\Om)$.
\end{theorem}
\noindent {\sc Proof:} We decompose the solution $\disp (v,p)$ of (\ref{eq3.37})
setting
\begin{equation} \label{eq3.39}
\left|
\begin{array}{l}
v = v_0 + g\\
p = p_0 + q,
\end{array}
\right.
\end{equation}
where $v_0$, $p_0$ is given by
\begin{equation} \label{eq3.40}
\left|
\begin{array}{l}
\disp v_0'' + L\,v_0 = 0 \ \mbox{ in } \ Q\\[5pt]\disp
v_0 = \left\{
\begin{array}{l}
0 \ \mbox{ on } \ \Sigma_1\\[5pt]
\disp \frac{1}{\sigma k^2(t)}\,({p_0})_{\nu} \ \mbox{ on } \ \Sigma_2\\[5pt]\disp
0 \ \mbox{ on } \ \Sigma \backslash \Sigma_0
\end{array}
\right.\\[5pt]\disp
v_0(0) = v_0'(0) = 0, \;\; y\in \Om,
\end{array}
\right.
\end{equation}
\begin{equation} \label{eq3.41}
\left|
\begin{array}{l}
\disp p_0'' + L^{\ast}p_0 = k^n(t) \left(v_0 - v_2\right) \ \mbox{ in } \ Q \\[5pt]\disp
p_0 = 0 \ \mbox{ on } \ \Sigma\\[5pt]\disp
p_0(T) = p_0'(T) = 0, \;\; y \in \Om,
\end{array}
\right.
\end{equation}
and $\disp\{g,q\}$ is given by
\begin{equation} \label{eq3.42}
\left|
\begin{array}{l}
g'' + L\,g = 0 \ \mbox{ in } \ Q\\[5pt]\disp
g = \left\{
\begin{array}{l}
w_1 \ \mbox{ on } \ \Sigma_1\\[5pt]
\disp \frac{1}{\sigma k^2(t)}\,q_\nu \ \mbox{ on } \ \Sigma_2\\[5pt]\disp
0 \ \mbox{on} \ \Sigma \backslash \Sigma_0
\end{array}
\right.\\[5pt]\disp
g(0) = g'(0) = 0, \;\;y \in \Om,
\end{array}
\right.
\end{equation}
\begin{equation} \label{eq3.43}
\left|
\begin{array}{l}
\disp q'' + L^{\ast}q = k^n(t) g \ \mbox{ in } \ Q \\[5pt]\disp
q = 0 \ \mbox{ on } \ \Sigma\\[5pt]\disp
q(T) = q'(T) = 0, \;\; y \in \Om.
\end{array}
\right.
\end{equation}
We set next
\begin{equation} \label{eq3.44}
\begin{array}{cccc}
A \ : & \! L^2(\Sigma_1) & \! \longrightarrow & \! H^{-1}(\Om) \times L^2(\Om) \\
& \! w_1 & \! \longmapsto & \! A\,w_1 = \big\{ g'(T;w_1) + \delta g(T;w_1),\; -g(T;w_1) \big\},
\end{array}
\end{equation}
which defines
$$A \in \mathcal{L}\left( L^2(\Sigma_1); \;H^{-1}(\Om) \times L^2(\Om)\right),$$
where  $\delta$ is  a positive constant .

Then (\ref{subj1}), using (\ref{eq3.39}) and (\ref{eq3.44}) can be rewritten as
\begin{equation} \label{subj2}
\begin{array}{l}
\disp Aw_1\in \{ -v_0(T)+\delta g(T)+B_{H^{-1}(\Om)}(v^1,\rho_1),\;-v_0(T)+B_{L^2(\Om)}(v^0,\rho_0)\}.
\end{array}
\end{equation}
We will show that  $Aw_1$ generates a dense subspace of $H^{-1}(\Om) \times
L^2(\Om)$. For this, let  $\{ f^0,f^1 \} \in H_{0}^{1}(\Om) \times L^2(\Om)$ and we
consider the following systems (''adjoint states"):
\begin{equation} \label{eq3.45}
\left|
\begin{array}{l}
\varphi'' + L^{\ast}\,\varphi =\disp k^n(t) \;\psi \ \mbox{ in } \ Q \\[5pt]\disp
\varphi = 0 \ \mbox{ on } \ \Sigma\\[5pt]\disp
\varphi(T) = f^0, \ \varphi'(T) = f^1, \;\; y \in \Om,
\end{array}
\right.
\end{equation}

\begin{equation} \label{eq3.46}
\left|
\begin{array}{l}
\psi'' + L\,\psi = 0 \ \mbox{ in } \ Q\\[5pt]\disp
\psi = \left\{
\begin{array}{l}
0 \ \mbox{ on } \ \Sigma_1\\[5pt]\disp
\disp \frac{1}{\sigma k^2(t)}\,\varphi_\nu \ \mbox{ on } \ \Sigma_2\\[5pt]\disp
0 \ \mbox{ on } \ \Sigma \backslash \Sigma_0
\end{array}
\right.\\[5pt]\disp
\psi(0) = \psi'(0) = 0, \;\; y \in \Om.
\end{array}
\right.
\end{equation}
Multiplying $(\ref{eq3.46})_1$ por $q$, $(\ref{eq3.45})_1$ by $g$, where $q$, $g$
solves (\ref{eq3.43}) and (\ref{eq3.42}), respectively, and integrating in $\disp Q=
\Om\times (0,T)$ we obtain
\begin{equation} \label{eq3.47}
\int_{0}^{T}\int_{\Om} k^n(t)g\,\psi\,dy\,dt =- \frac{1}{\sigma}\int_{\Sigma_2} \frac{1}{k^4(t)}\,q_\nu\, \varphi_\nu d\Sigma,
\end{equation}
and
\begin{equation} \label{eq3.49}
\disp \langle g'(T),f^0 \rangle_{H^{-1}(\Om),H_{0}^{1}(\Om)} + \delta \langle g(T), f^0 \rangle_{L^2(\Om),H_{0}^1(\Om)}  - \big( g(T),f^1 \big) =-\int_{\Sigma_1}\frac{1}{k^2(t)}\,\varphi_\nu\,w_1\,d\Sigma.
\end{equation}
Consider the first hand side of (\ref{eq3.49}) as the inner product of $\disp \{g'(T)+
\delta g(T),-g(T)\}$ with $\{ f^0,f^1 \}$ in $ H^{-1}(\Om)  \times L^2(\Om) $ and $
H_{0}^{1}(\Om) \times L^2(\Om)$.

We then obtain from (\ref{eq3.49})
\begin{equation*}
\Big\langle \big\langle A\,w_1 , f \big\rangle \Big\rangle = - \int_{\Sigma_1}\frac{1}{k^2(t)}\,\varphi_\nu\,w_1\,d\Sigma,
\end{equation*}
where $\Big\langle \big\langle . , . \big\rangle \Big\rangle$ represent the duality
pairing between $ H^{-1}(\Om) \times L^2(\Om) $ and $ H_{0}^{1}(\Om) \times
L^2(\Om) $. Therefore if
$$\langle g'(T),f^0 \rangle_{H^{-1}(\Om),H_{0}^{1}(\Om)} + \delta \langle g(T), f^0 \rangle_{L^2(\Om),H_{0}^1(\Om)} - \big( g(T),f^1 \big) = 0,$$
for all $w_1 \in L^2(\Sigma_1)$, then
\begin{equation} \label{eq3.50}
\disp \varphi_\nu= 0 \ \ \mbox{ on } \ \ \Sigma_1.
\end{equation}
Then in case (\ref{decomp 2.A}),
\begin{equation} \label{eq3.51}
\psi = 0 \ \ \mbox{ on } \ \ \Sigma, \;\; \mbox{ so that } \psi\equiv 0.
\end{equation}
Therefore
\begin{equation} \label{eq3.54}
\begin{array}{l}
\disp \varphi'' + L^{\ast}\,\varphi = 0, \;\; \varphi = 0 \mbox{ on }  \Sigma,
\end{array}
\end{equation}
and satisfies (\ref{eq3.50}). Therefore, according to  Holmgren's Uniqueness
Theorem (cf. \cite{LH})   and for the explicit use of it made here (cf. \cite {Mi})
and if (\ref{hT}) holds, $\disp \varphi \equiv 0$, so that $\disp f^0=0, f^1=0$. This
ends the proof.\Fin

\section{Optimality system for the leader}\label{sec5}
Thanks to the results obtained in the Section \ref{sec3}, we can take, for the
each $\disp w_1$, the Nash equilibrium $\disp w_2$ associated to solution
$\disp v$ of (\ref{eq1.14}). Next, for completeness, we will show that there exists
a leader control $\disp w_1$ solution of the following problem:
\begin{equation} \label{eq3.7cil.1}
\disp \inf_{w_1\in \mathcal{U}_{ad}} J(w_1)
\end{equation}
where $\disp \mathcal{U}_{ad}$ is the set of admissible controls
\begin{equation}\label{admcon}
\disp \mathcal{U}_{ad}=\{w_1\in L^2({\Sigma_1}); \; v \mbox{ solution of } (\ref{eq1.14}) \mbox{ satisfies } (\ref{subj1})\}.
\end{equation}

 For this, we will a duality argument due to Fenchel and Rockfellar \cite{R} (see also,
\cite{Bre}, \cite{EK}).
\begin{theorem} \label{teor3.6} Assume the hypotheses $\disp (H1) - (H4)$, (\ref{decomp 2.A}) and (\ref{hT}) are satisfied. Then for $\{f^0,f^1\}$ in $\disp H_0^1(\Om) \times L^2(\Om)$ we uniquely define $\{\varphi, \psi, v, p \}$ by
\begin{equation} \label{eq3.139}
\left|
\begin{array}{l}
\disp \varphi'' + L^* \varphi = k^n(t)\psi \ \ \text{in} \ \ Q \\[3pt]\disp
\disp \psi'' + L \psi = 0 \ \ \text{in} \ \ Q \\[3pt]\disp
\disp v'' + Lv = 0 \ \ \text{in} \ \ Q \\[3pt]\disp
\disp p'' + L^*p  = k^n(t)(v - v_2) \ \ \text{in} \ \ Q \\[3pt]\disp
\disp \varphi = 0 \ \ \ \text{on} \ \ \ \Sigma \\[3pt]\disp
\disp \psi =
\left\{
\begin{array}{l}
\disp 0 \ \ \text{on} \ \ \Sigma_1 \\[3pt]\disp
\disp \frac{1}{\sigma k^2(t)}\;\varphi_\nu \ \ \text{on} \ \ \Sigma_2 \\[3pt]\disp
\disp 0 \ \ \ \text{on} \ \ \ \Sigma \backslash \Sigma_0\\[3pt]\disp
\end{array}
\right.\\
\disp v =
\left\{
\begin{array}{l}
\disp - \frac {1}{k^2(t)}\;\varphi_\nu \ \ \text{on} \ \ \Sigma_1\\[3pt]\disp
\disp \frac{1}{\sigma k^2(t)}\;p_\nu \ \ \text{on} \ \ \Sigma_2\\[3pt]\disp
\disp 0 \ \ \ \text{on} \ \ \ \Sigma \backslash \Sigma_0\\[3pt]\disp
\end{array}
\right.\\
\disp p = 0 \ \ \ \text{on} \ \ \ \Sigma \\[3pt]\disp
\disp \varphi(.,T) = f^0,\, \varphi'(.,T) = f^1 \ \ \text{in} \ \  \Om \\[3pt]\disp
\disp v(0) = v'(0) = 0 \ \ \text{in} \ \   \Om \\[3pt]\disp
\disp p(T) = p'(T) = 0 \ \ \text{in} \ \  \Om
\end{array}
\right.
\end{equation}
We uniquely define  $\{f^0,f^1\}$, as the solution of the variational inequality
\begin{equation} \label{eq3.140}
\begin{array}{l}
\disp \big\langle v'(T,f) - v^1, \widehat{f}^0 - f^0\big\rangle_{H^{-1}(\Om), H_{0}^{1}(\Om)} - \big(v(T,f) - v^0, \widehat{f}^1 - f^1\big)  \\[10pt]
\disp + \rho_1\big(||\widehat{f}^0|| - ||f^0||\big) + \rho_0\big(|\widehat{f}^1| - |f^1|\big) \geq 0,\,\forall\, \widehat{f} \in H_{0}^{1}(\Om) \times L^2(\Om).
\end{array}
\end{equation}
Then the optimal leader is given by
\begin{equation*}
w_1 = -\frac{1}{k^2(t)}\;\varphi_\nu \ \ \text{on} \ \ \Sigma_1,
\end{equation*}
where $\varphi$ corresponds to the solution of (\ref{eq3.139}). 
\end{theorem}
\noindent {\sc Proof:} We introduce two convex proper functions as follows, firstly
\begin{equation}\label{eq3.119}
\begin{array}{l}
\disp F_1 : L^2(\Sigma_1) \longrightarrow \re \cup \{\infty\},\\[5pt]
\disp F_1(w_1) = \frac{1}{2} \int_{\Sigma_1}w_{1}^{2}\,d\Sigma
\end{array}
\end{equation}
the second one
\begin{equation*}
F_2 : H^{-1}(\Omega) \times L^2(\Omega) \longrightarrow \re \cup \{\infty\}
\end{equation*}
given by
\begin{align} \label{eq3.120}
\nonumber F_2(Aw_1) &= F_2\big(\{g'(T,w_1) + \delta g(T,w_1),-g(T,w_1)\}\big) =\\
& = \left\{
\begin{array}{l}
0, \text{ if }
\left\{
\begin{array}{l}
g'(T) + \delta g(T) \in v^1 - v_0'(T) + \delta g(T) + \rho_1B_{H^{-1}(\Om)}\\
-g(T) \in -v^0 + v_0(T,w_1) - \rho_0B_{L^2(\Om)}
\end{array}
\right.\\
+ \infty, \text{ otherwise}.
\end{array}
\right.
\end{align}
With these notations problems (\ref{inf1}), (\ref{subj1}) becomes equivalent to
\begin{equation} \label{eq3.122}
\begin{array}{l}
\disp \inf_{w_1 \in L^2(\Sigma_1)}\big[F_1(w_1) + F_2(Aw_1)\big]
\end{array}
\end{equation}
provided we prove that the range of $\disp A$ is dense in $\disp H^{-1}(\Om)
\times L^2(\Om)$, under conditions (\ref{hT}).

By the Duality Theorem of Fenchel and Rockfellar \cite{R}( see also \cite{Bre},
\cite{EK}), we have
\begin{equation} \label{eq3.124}
\begin{array}{l}
\inf_{w_1 \in L^2(\Sigma_1)}[F_1(w_1) + F_2(Aw_1)]\\[5pt]\disp  = -\inf_{(\widehat{f}^0,\widehat{f}^1) \in H_{0}^{1}(\Om) \times L^2(\Om)} [F_{1}^{*}\big(A^*\{\widehat{f}^0,\widehat{f}^1\}\big) + F_{2}^{*}\{-\widehat{f}^0, -\widehat{f}^1\}],
\end{array}
\end{equation}
where $\disp F_i^*$ is the conjugate function of $\disp F_i  (i=1,2)$ and $\disp
A^*$ the adjoint of $\disp A$.

We have
\begin{equation} \label{eq3.121}
\begin{array}{ccccc}
A^* \ : & \! H_{0}^{1}(\Omega) \times L^2(\Omega) & \! \longrightarrow & \! L^2(\Sigma_1) \\
& \! (f^0,f^1) & \! \longmapsto & \! A^*f = & \! -\dfrac{1}{k^2}\,\varphi_\nu,
\end{array}
\end{equation}
where $\varphi$ is give in (\ref{eq3.45}).
We see easily that
\begin{equation} \label{eq3.125}
F_{1}^{*}(w_1) = F_1(w_1)
\end{equation}
and
\begin{equation}\label{eq3.125.2}
\begin{array}{l}
\disp  F_{2}^{*}(\{\widehat{f}^0,\widehat{f}^1 \})  = \big( v_0(T) - v^0 ,\widehat{f}^1\big) + \langle v^1 - v_0'(T) + \delta g(T), \widehat{f}^0\rangle_{H^{-1}(\Omega) \times H_{0}^{1}(\Omega)}\\ [10pt]
\disp + \rho_1||\widehat{f}^0|| + \rho_0|\widehat{f}^1|
\end{array}
\end{equation}
Therefore the (opposite of) right hand side of (\ref{eq3.124}) is given by
\begin{align} \label{eq3.127}
\disp
& - \inf_{\widehat{f} \in H_{0}^{1} \times L^2(\Omega)} \bigg\{\frac{1}{2}\int_{\Sigma_1}\left(\frac{1}{k^2(t)}\right)^2 \widehat{\varphi}_\nu^2d\;\Sigma  + \big( v^0 - v_0(T) ,\widehat{f}^1\big) \\
\nonumber & - \langle v^1 - v_0'(T) + \delta g(T), \widehat{f}^0\rangle_{H^{-1}(\Omega) \times H_{0}^{1}(\Omega)} + \rho_1||\widehat{f}^0|| + \rho_0|\widehat{f}^1|\bigg\}
\end{align}
This is the dual problem of (\ref{inf1}), (\ref{subj1}).

We have now two ways to derive the optimality system for the leader control,
starting from the primal or from the dual problem. \Fin

\end{document}